\def\F{\mathbb F}
  
    \def\N{\mathbb N}
    
  \def\a{\alpha}
  \def\b{\beta}

\def\k{\kappa}

  \def\e{\epsilon}

\def\inn{\mathrm{Inn}}
\def\aut{\mathrm{Aut}}
  
  \def\la{\langle}
  \def\ra{\rangle}

    \documentclass[a4paper]{article}
            \usepackage{amsfonts,amssymb}
\usepackage{amsmath}

\newtheorem{thm}{Theorem}
  \newtheorem{prop}[thm]{Proposition}
  \newtheorem{lem}[thm]{Lemma}
  \newtheorem{cor}[thm]{Corollary}

\begin{document}

  \title{Verbal width in anabelian groups}
  \author{Nikolay Nikolov }
    \date{}
  \maketitle

 \begin{abstract}
The class $\mathcal A$ of anabelian groups is defined as the collection of finite groups without abelian composition factors. We prove that the commutator word $[x_1,x_2]$ and the power word $x_1^p$ have bounded width in $\mathcal A$ when $p$ is an odd integer. By contrast the word $x^{30}$ does not have bounded width in $\mathcal A$. On the other hand any given word $w$ has bounded width for those groups $G \in \mathcal A$ whose composition factors are sufficiently large as a function of $w$. In the course of the proof we establish that sufficiently large almost simple groups cannot satisfy $w$ as a coset identity.
 \end{abstract}
\section{Introduction} The study of word values and verbal width of groups has received considerable attention recently, see for example \cite{dan} and the survey \cite{shalev}. We say that a word has width at most $m$ in a group $G$ if every element of the verbal subgroup $w(G)$ is a product of at most $m$ values of $w$. It turns out that finite nonabelian simple groups have bounded word width: in fact by \cite{LS} any given word has width at most 2 in every sufficiently large finite simple group. For general finite groups it was proved in \cite{DanNN} and \cite{NS} that the commutator word $[x_1,x_2]$ and the power word $x_1^q$ have bounded width in the class of $d$-generated finite groups.
\footnotetext{
  The author is grateful for the support of an EPSRC grant.
  }
We could ask if the above statements hold true without the condition on the number of generators of the finite groups. It is not too hard to see that  $[x_1,x_2]$ and $x_1^p$  have unbounded width even in the class of finite $p$-groups. 

More generally we have the following result. A word $w$ in the free group $F$ is called \emph{silly} if $w=1$ or if $w(F)=F$.
\begin{thm} \label{silly} Let $w$ be word which is not silly. Then $w$ has unbounded width in the class of all finite groups.
\end{thm}

However the situation is different when one considers finite groups without abelian composition factors. We will call such groups \emph{anabelian}. This class is a natural place to try to apply and extend the current knowledge about word values in finite simple groups. Our next result is a relatively straightforward consequence of Proposition 11.1 from \cite{DanNN}.

\begin{thm}\label{com} There is a constant $D \in \N$ such that every element of a finite anabelian group is a product of $D$ commutators.
\end{thm}
The affirmative resolution of Ore conjecture \cite{ore} proved that every element in a finite nonabelian simple group is a commutator. We expect that the constant $D$ in the above Theorem can be taken to be equal to 1. Some explicit upper bounds for $D$ can be found in \cite{levy}.

The corresponding question for powers has a more complicated answer.
\footnotetext{2010 {\it Mathematics Subject Classification:} 20D40,
  20D06 }
\begin{thm} \label{p} Given an odd integer $p$ there is an integer $l=l(p) \in \N$ such that every element of an anabelian group $G$ is a product of at most $l$ $p$-th powers.
\end{thm}

\begin{thm} \label{60} Given integer $m \in \N$ there is an anabelian group $G$ and an element $g \in G=G^{30}$ which is not a product of less than $m$ $30$-th powers. \end{thm}
The integer $30$ in the above theorem can be replaced with the exponent of any nonabelian finite simple group.

In fact the presence of small composition factors is the only reason verbal width is unbounded in anabelian groups.

\begin{thm} \label{w} Given a group word $w$ there are integers $f=f(w)$ and $c=c(w) \in \N$ such that $w$ has width at most $f$ in any anabelian group whose composition factors have size at least $c$.
\end{thm}  
We conjecture that the number $f$ in Theorem \ref{w} does not depend on $w$ and can be taken to be a small constant less than 10. For example the proof of Proposition \ref{key} can be modified to show that in the special case when $G$ is an iterated wreath product of semisimple groups with sufficiently large simple factors then $f$ can be taken to be 5.

The proof of Theorem \ref{w} follows by induction on $|G|$ from the following general result. 

\begin{prop} \label{key} Given a group word $w=w(x_1, \ldots, x_d)$ there are integers $f=f(w)$ and $c=c(w) \in \N$ with the following property. Let $G$ be a finite group with a normal semisimple subgroup $N=S^{(k)}$ and let $\mathbf a_1, \ldots, \mathbf a_f \in G^{(d)}$ be arbitrary $d$-tuples in $G$. The map $\psi : N^{df} \rightarrow N$ defined by
\[ \psi (\mathbf z_1, \ldots, \mathbf z_f)= \left( \prod_{i=1}^f w(\mathbf a_i \mathbf z_i)  \right) \left( \prod_{i=1}^f w(\mathbf a_i)  \right)^{-1}\]
is surjective provided $|S|>c$.
\end{prop}

Note that we do not require that $G$ is anabelian in the above Proposition. Its proof depends among other things on the following result which may have independent interest.

\begin{prop} \label{coset} Let $w(x_1, \ldots, x_d)$ be a group word. There is an integer $c_0=c_0(w)$ with the folllowing property: 

Let $G$ be a finite group with semisimple normal subgroup $N=S^{(m)}$ and let $g_1, \ldots, g_d \in G$. Assume that $|S|>c_0$. There exist $a_1, \ldots, a_d \in N$ such that $w(g_1a_1, \ldots, g_da_d)$ does not centralize any simple factor of $N$. 
\end{prop}

This proposition implies that a given coset identity cannot hold in any large enough almost simple group.
\begin{cor} \label{cos} Given a group word $w(x_1, \ldots, x_d)$ let $c_0=c_0(w)$ be the integer proved by Propostion \ref{coset} above. Then if $G$ is an almost simple group with socle $S$ of size greater than $c_0$ and $g_1, \ldots, g_d \in G$ then $w(g_1S, \ldots , g_dS) \not = \{1\}$.
\end{cor}
    
	\section{Proofs}
	Generally the proofs of Theorems \ref{silly}, \ref{com}, \ref{p} and \ref{60} follow quickly from results in \cite{DanNN}. On the other hand Propositions \ref{key} and \ref{coset} require some new technical results on product decompositions and equations in almost simple groups.

\bigskip

\subsection{Proofs of Theorems \ref{silly} and Theorem \ref{60}}

Let us prove Theorem \ref{60} first.
Let $L$ be an anabelian group such that $L^{30}=L$ and which cannot be generated by $m$ elements, for example we can take $L= A_6^{(N)}$ where $N> |A_6|^m$. For each $m$-tuple $\mathbf y =(y_1 ,\ldots, y_m) \in L^{(m)}$ let $\Omega_\mathbf{y}:=L/\la y_1, \ldots, y_m\ra$ be set of right cosets of the proper subgroup $\la y_1, \ldots, y_m\ra$ in $L$. Let $N_\mathbf y = A_5^{\Omega_{\mathbf y}}$ and define action of $L$ on $N_{\mathbf y}$ by permuting the factors according to the transitive action of $L$ on the coset space $\Omega_{\mathbf y}$. 
Let $K= \bigoplus_{ \mathbf y \in L^{(m)}} N_{\mathbf y}$  and take $G:= K \rtimes L$. Since $K$ is a product of the minimal normal subgroups $N_{\mathbf y}$ and $L=L^{30}$ we have that $G=G^{30}$. On the other hand let $\k \in K$ be any element with nontrivial projections on every factor $A_5$ in every $N_{\mathbf y}$.

Suppose that $\kappa= g_1^{30} \cdots g_m^{30}$ for some $g_i \in G$. Write $g_i= k_i l_i$ with $k_i \in K, l_i \in L$ and let $\mathbf a= (l_1, \ldots, l_m)$, $\omega= \la l_1, \ldots, l_m \ra \in \Omega_{\mathbf a}$. Note that the direct factor $A \simeq A_5$ of $K$ labelled by $\omega$ in $N_{\mathbf a}= A_5^{\Omega_{\mathbf a}}$ is centralized by all $l_i$.

Let us denote by $u(\omega)$ the projection of any element $u \in K$ onto the direct factor $A$. We have $\kappa (\omega)= k_1(\omega)^{30} \cdots k_1(\omega)^{30}$. Since $A^{30}=1$ we conclude that $\kappa (\omega)=1$, contradiction. $\square$ \bigskip

The proof of Theorem \ref{silly} is similar. Let $w=w(x_1, \ldots, x_k)$ be a word in the free group $F_k$ with basis $x_1, \ldots, x_k$. For $i=1, \ldots, k$ let $c_i$ be the exponent sum of the occurences of the letter $x_i$ in $w$. Note that since $w$ is not silly there is a positive integer $d>1$ which divides each $c_i$. Let $C$ be the cyclic group of order $d$, we have $w(C)=1$. There is a finite simple group $B$ such that $w(B)\not =1$ and hence $w(B)=B$. 

Let $n \in \mathbb  N$. We will exhibit a finite group $G$ such that $w$ has width at least $n+1$ in $G$. Choose an integer $m$ such that $L:=B^{(m)}$ cannot be generated by $nk$ of its elements. For $\mathbf y \in L^{(nk)}$ let $\Omega_{\mathbf y} = L/ \la y_1 , \ldots, y_{nk}\ra$ and let $N_{\mathbf y}= C^{\Omega_{\mathbf y}}$. As before we take $K= \bigoplus_{ \mathbf y \in L^{(nk)}} N_{\mathbf y}$ and define $G= K \rtimes L$. Since $L=w(L)$ and $K$ is a product of minimal normal subgroups of $G$ with trivial centralizers in $G$, we deduce that $G=w(G)$. 

Let $\kappa \in K$ be an element with nontrival projections onto every factor $C$ of each $N_{\mathbf y}$. We claim that $\kappa$ is not a product of $n$ values of $w$. Suppose that
$\kappa = \prod_{i=1}^n w(\mathbf z_i)$ where $\mathbf z_i= (z_{i,j})_{j=1}^k$ with $z_{i,j} \in G$. Write each $z_{i,j}= k_{i,j}l_{i,j}$ with $k_{i,j} \in K, l_{i,j} \in L$ and let $\mathbf l$ be the $nk$-tuple of the elements $\{l_{i,j} \ | \ i=1, \ldots, n, \ j=1, \ldots, k \}$. Since $w(C)=1$ the argument at the end of the proof of Theorem \ref{60} gives that  $\kappa$ must have at least one trivial coordinate in the subgroup $N_{\mathbf l}$ corresponding to $\mathbf l$. Therefore $w$ has width at least $n+1$ in $G$. $\square$

\subsection{Proof of Theorem \ref{com}.}
Our proof of Theorem \ref{com} is a straighforward consequence of Proposition \ref{ss} below.
	
First we need to intruduce some notation: For automorphisms $\alpha, \beta$ of a group $G$ and $x,y \in G$ define $T_{\a,\b}(x,y)=x^{-1}y^{-1}x^\a y^\b$. For two $n$-tuples $\mathbf a =(\a_1, \ldots, \a_n), \mathbf b =(\b_1, \ldots ,\b_n)$ in $\aut(G)^{(n)}$ let $T_{\mathbf{a,b}}(\mathbf{x,y})= \prod_{i=1}^n T_{\a_i,\b_i}(x_i,y_i)$ for $x_i,y_i \in G$. Define \[ T_{\mathbf{a,b}}(G,G)= \{ T_{\mathbf{\mathbf{a,b}}}(\mathbf{x,y})\ | \ \forall \mathbf{x,y} \in G^{(n)}\} \]
Finally define $[\a,G]:= \{ [\a, g] \ | \ g \in G\}$.

We will use the following result from \cite{DanNN}, Proposition 11.1.
\begin{prop}\label{ss} There is a constant $D \in \N$ with the following property. Let $N$ be a semisimple group, i.e. a product of finite simple groups and let $\mathbf{a,b} \in \aut (N)^{(D)}$. Then $T_{\mathbf{a,b}}(N,N)=N$.
\end{prop}

We can now prove Theorem \ref{com} by induction on the size of the anabelian group $G$. When $G$ is simple the Theorem follows from the validiy of Ore's conjecture which states that every element of $G$ is a commutator.

In general consider a minimal normal subgroup $N$ of $G$ and an element $g \in G$. By the induction hypothesis we may assume that there are elements
$u_i,v_i \in G, \ i=1, \ldots D$ such that $g= \prod_{i=1}^D[u_i,v_i] \k$ for some $\k \in N$. Now we search for elements $x_i,y_i \in N$ such that
$g= \prod_{i=1}^D[x_iu_i,y_iv_i]$. Lemma 4.6 of \cite{DanNN} gives that
\[ (\prod_{i=1}^D[x_iu_i,y_iv_i]) (\prod_{i=1}^D[u_i,v_i])^{-1}= \prod_{i=1}^D T_{a_i,b_i}(x_i^{\sigma_i},y_i^{\rho_i}) \] for certain elements $a_i,b_i,\sigma_i,\rho_i \in G$ which depend only on $u_1, \ldots, u_D, v_1, \ldots ,v_D$. By Proposition \ref{ss} the equation 
$\prod_{i=1}^D T_{a_i,b_i}(x_i^{\sigma_i},y_i^{\rho_i})= \k$ has a solution in $x_i,y_i \in N$ and the induction step is complete. $\square$

\subsection{Proof of Theorem \ref{p}.}

We will need the following result which may have independent interest.

\begin{prop} \label{p'}Let $G \leq \mathrm{Aut}(S)$ be an almost simple group with simple socle $S$. For any $a\in G$ and an odd integer $n$ there is some $g \in S$ such that $(ag)^n \not =1$.
\end{prop}

\textbf{Proof:}

We will use \cite{GLS} as a reference for information on the finite simple groups and their automorphisms. In particular for a finite simple group of Lie type $S$ defined over  a field $\F$ we shall denote by $D_S, \Phi_S$ the subgroups of $\aut(S)$ of diagonal, field, and a subset $\Gamma_S$ of graph automorphisms. ($\Gamma_S$ is a subgroup of $\aut(S)$ unless $S$ is a Suzuki or Ree group). 
We will omit the subscript $S$ and write $D, \Phi$ and $\Gamma$ when the simple group $S$ is clear from the context. We have $\aut(S)= \mathrm{Inn}(S) D\Phi \Gamma$.
For a prime power $q$ we denote by $[q] \in \Phi$ the field automorphism  of $S$ induced by the automorphism $x \mapsto x^q$  of $\F$.

Suppose that the statement of the proposition is false, i.e there is some $a \in G$ such that $(ag)^n=1$ for all $g \in S$. In the first place note that $a^n=1$ and since every finite simple group has an elements of even order we conclude that $a \not \in S$. Therefore $a$ has order odd order at least 3 in $G/S$.

 This implies that $S$ cannot be a sporadic or an alternating simple group, because in this case $\mathrm{Out}(S)$ has exponent 2.
So $S$ must be simple group of Lie type. By replacing $a$ with $as$ for appropriate $s \in S$ we may assume that $a$ acts on $S$ as an element of the subgroup $D \Phi \Gamma \leq \aut(S)$ generated by the diagonal, field and graph automorphisms. 

Suppose first that $S$ is not of type $^2B_2$. We claim that $S$ contains a subgroup $S_0$ isomorphic to $(P)SL(2,q)$ such that $S_0^a=S_0$. The subgroup $\aut_0(S):= \mathrm{Inn}(S) D\Phi$ has index $2$ in $\aut (S)$ unless $S= D_4(q)$ and since the order of $a$ in $G/S \leq \mathrm{Out}(S)$ is odd we conclude that in this case $a$ acts on $S$ as an element of $\aut_0(S)$ i.e. without a graph automorphism component. Therefore $a$ preserves all root subgroups of $S$ and in particular $a$ stabilizes a quasi-simple subgroup $S_0$ of $S$ of type $A_1$. If $S$ has type $D_4$ it is evident that $D \Phi \Gamma$ preserves the central root subgroup of the Dynkin diagram and again $a$ preserves a copy $S_0$ of $(P)SL_2(q)$ inside $S$.

In the equation $(aa)^n=1$ we restrict $g \in S_0$ and by considering $\langle S_0,a \ra$ acting on itself by conjugation we are reduced to the case when $S=PSL_2(q)$.
Now the subgroup $\inn (S)\Phi_S$ generated by all inner and field automorphisms is a normal subgroup of index at most $2$ in $\aut(S)$. Using again that $a$ has odd order in $G/S \leq \mathrm{Out}(S)$ we conclude that $a \in \inn(S) \Phi$. Again, replacing $a$ by appropriate $as$ with $s \in S$ we may assume that $a$ is a field automorphism of $S$. Now consider the equation $(ag)^n=1$ with $g \in PSL_2(q_0)$ where $\mathbb F_{q_0}$ is the ground field of $\mathbb F_q$. Since $a$ centralizes $g$ we deduce that $g^n=1$, and therefore $PSL_2(q_0)$ has exponent $n$. This is a contradiction since $|PSL_2(q_0)|$ is always even. 

Finally, if $S= \ ^2B_2$ then $S$ has no outer diagonal automorphism and we may assume that $a \in \Phi$. If $g \in \ ^2B_2(2) < S$ then $a$ centralizes $g$ and therefore $(ag)^n=a^n=1$ implies $g^n=1$ for all $g \in ^2B_2(2)$. However $^2B_2(2)$ is a Frobenius groups of size 20 and does not have odd exponent.
$\square$ 
\bigskip

The following are Propositions 10.1 and 9.1 of \cite{DanNN}.
\begin{prop}\label{10.1} For every integer $q$ there is an integer $C=C(q)$ with the following property. 
Let $N$ be a semisimple normal subgroup of a group $G$ and let $h_1, \ldots h_m \in G$. Assume that $m>C$ and each finite simple factor of $N$ has size at least $C$. The mapping $\psi : N^{(m)} \rightarrow N$ defined by 
\[ \prod_{i=1}^m (x_ih_i)^q= \psi (x_1,\ldots,x_m) \prod_{i=1}^m h_i^q\]
is surjective.
\end{prop}

\begin{prop}\label{9.1} Let $S$ be a finite simple group and let $T=S^{(n)}$. Let $f_1, \ldots, f_m \in \mathrm{Aut}(T)$ and let $c(f_i)$ denote the number of cycles of $f_i$ on the $n$ simple factors of $T$. Assume that $\langle f_1, \ldots, f_m\rangle $ induces a transitive permutation group on the simple factors of $T$. There is a constant $D$ such that if $\sum_{i=1}^m c(f_i) \leq n(m -2) -2D$ the map $\phi : T^{(n)} \rightarrow T$
defined by
\[ \phi(x_1, \ldots, x_m)= \prod_{i=1}^m [f_i,x_i], \quad x_i \in T \]
is surjective.
\end{prop}

We will also need the following.

\begin{prop} \label{fix} Let $N=S^{(k)}$ be a semisimple normal subgroup of a group $G$, let $h \in G$ and let $p$ be an odd integer. There is an element $b \in N$ with the following property: $(hb)^p$ does not centralize any simple factor of $N$. 
\end{prop}
\medskip

\textbf{Proof.} Without loss of generality we may assume that $h$ acts by conjugation on the $k$ simple factors of $N=S_1 \times \cdots S_k$ as the $k$ cycle $S_i^h=S_{i+1}$ for $i=1, \ldots, k$ (where $S_{k+1}=S_1$). If $k$ does not divide $p$ then $h^p$ does not stabilize any factor of $N$ and we take $b=1$. So suppose $p=kp_1$ where $p_1 \in \mathbb N$.  Let $h_1=h^k$ and note that $h_1$ stabilizes each $S_i$. Let $a_i \in \mathrm{Aut}(S_i)$ be the automorphism of $S_i$ induced by conjugation by $h_1$. 

It is sufficient to find $b \in N$ such that $(hb)^p$ does not centralize $S_1$: Since $hb$ permutes the factors $S_i$ cyclically and commutes with $(hb)^p$ it will then follow that $(bh)^p$ does not centralize any of the $S_i$.

Let $u \in S_1$ be the element such that $(a_1u)^{p_1} \not =1$ provided by Proposition \ref{p'}. Let $b=(u,1,1, \ldots, 1) \in S^{(k)} =N$ and note that 

\[ (hb)^k= h^k b^{h^{k-1}} \cdots b^h b\] acts on $S_1$ as the automorphism $h^k u=a_1u$
Therefore $(hb)^n= ((hb)^k)^{p_1}$ acts on $S_1$ as $(a_1u)^{p_1}$ which is a nontrivial automorphism of $S_1$ by the choice of $u$.
$\square$

We are going to prove Theorem \ref{p} by induction on the size of the anabelian group $G$. When $G$ is simple the theorem follows from 
the results in \cite{SW} and \cite{MZ} as well as the stronger result in \cite{LS}. Let $N$ be a minimal normal subgroup of $G$ isomorphic to $S^{(n)}$ where $S$ is a simple group. Let $C$ and $D$ be the number provided by Propositions \ref{10.1} and \ref{9.1} and take $l>4+C+4D$. 

Let $g \in G$ assume that we have found elements $h_1, \ldots, h_l \in G$ such that $g =  h_1^p \cdots h_l^p \kappa$ for some $\k \in N$. Suppose first that $|S|>C$. Since $l>C$ Proposition \ref{10.1} now implies that there exist $x_i \in N$ such that $g = (x_1h_1)^p \cdots (x_lh_l)^p$ and we are done.

Suppose now that $|S| <C$. According to Proposition $\ref{fix}$ we may replace $h_i$ with $h_ib_i$ for appropriate $b_i \in N$ and assume that
$h_i^p$ do not centralize any simple factor of $N$. We are now searching for elements $y_i \in N$ such that $g=(h_1^p)^{y_1} \cdots (h_l^p)^{y_l}$, i.e. \[ \kappa= (h_1^p \cdots h_l^p)^{-1} (h_1^p)^{y_1} \cdots (h_l^p)^{y_l}  \]
This equation is equivalent to 
\[ \kappa= \prod_{i=1}^l [f_i, z_i]  \quad \quad (\star) \]
where $f_i=(h_i^p)^{\sigma_i}, z_i=x_i^{\tau_i}$ for specific elements $\rho_i, \tau_i \in \la h_1^p, \ldots, h_l^p \ra$. Note that each the elements $f_i$ cannot centralize any simple factor of $N$.

We proceed to solve $(\star )$ independently in $\prod_{U \in \Omega} U  \leq N$ for each orbit $\Omega$ of $\la f_1, \ldots, f_l \ra$ on the simple factors of $N$. Therefore we may assume without loss of generality that the action of $\la f_1, \ldots, f_l \ra$ on the simple factors is transitive.

Let $n_i$ be the number of fixed points of $f_i$ on the simple factors of $N$ and $c_i=c(f_i)$ be the number of cycles. We have $c_i  \leq (n + f_i)/2$.
Suppose $\sum_{i=1}^l n_i < C$. In that case \[ \sum_{i=1}^l c_i  \leq \frac{1}{2}\sum_{i=1}^l(n + f_i) < \frac{nl +C}{2} \leq n(l-2) -2D,\]
where the last inequality holds because $l>4+C+4D$. So by Proposition \ref{9.1} we can solve $(\star )$ in $z_i \in N$ and we are done.

Finally, suppose that $\sum_{i=1}^l n_i > C$.
We can write the action of $f_i$ on $\mathbf a= (a_j) \in S^{(n)}=N$ as $\mathbf a^{f_i}= \mathbf b$ where $b_j=a_{\pi_i (j)}^{u_{i,j}}$ for a permutation $\pi_i\in \mathrm{Sym}(n)$ and automorphisms $u_{i,j} \in \aut(S)$.
Writing $\kappa= (k_j) \in N$ and $z_i= (z(i))_j$ the equation $(\star )$ becomes a system $\mathcal E$ of $n$ equations in $z(i)_j \in S$:

\[ E_j: \quad \quad k_j= \prod_{i=1}^l (z(i)_{\pi_i(j)})^{-u_{i,j}} \cdot z(i)_j
\]
Note that by our assumptions the permutations $\pi_1, \ldots , \pi_l$ generate a transitive subgroup on $\{1,\ldots, n\}$. Also each variable $z(i)_j$ appears exactly twice in all equations of $\mathcal E$.
We proceed to solve $\mathcal E$ by successively eliminating some of the elements $z(i)_j$ to reduce it to a single equation in $S$:

 At each step if we have more than one equation remaining we find a symbol $z(i)_j$ which appears in two different equations, say $E_r$ and $E_s$. We solve $E_r$ for $z(i)_j$ and substitute this value in $E_s$. The transitivity of $\pi_1, \ldots , \pi_l$ implies that this process will continue until we are left with a single equation $E_t: \ k_t=W$. Note that we have not eleimitated the letters $z(i)_j$ for those $(i,j)$ such that $\pi_i(j)=j$. Denoting $A:= \{ (i,j) \ | \ \pi_i(j)=j \}$ we find that the final equation has the form
\begin{equation} \label{finaleq} E_t: \quad k_t=\prod_{(i,j) \in A} X_{i,j} [u_{i,j}, z(i)_j] \cdot Y_{i,j} \end{equation} for some ordering of the set $A$ and some expressions $X_{i,j}$ and $Y_{i,j}$ which don't involve $z(i)_j$ with $(i,j) \in A$. Now choose at random and fix the values for the remaining variables in the expressions $X_{i,j},Y_{i,j}$ In this way $X_{i,j},Y_{i,j}$ become equal to some constants $d_{i,j}, e_{i,j} \in S$.
Let $T_{i_j}=[u_{i,j},S] := \{ [u_{i,j},s] \ | \ s \in S \}$. By assumption $u_{i,j} \in \aut(S)$ is not the identity and hence $|T_{i,j}|>1$. The existence of solutions $z(i)_j \in S$ to the final equation $E_t$ now follows from
\[S= \prod_{(i,j) \in A} d_{i,j} T_{i,j}e_{i,j}, \] which in turn is a direct consequence of the fact that $|A|>C >|S|$ together with the following 
\begin{prop} Let $S$ be a finite simple group. Let $v_i \in \aut(S) \backslash \{1\}$ and $r_i,r'_i \in S$ for $i=1,\ldots, m$. Assuming that $m > |S|$ we have
\[ S= \prod_{i=1}^m r_i [v_i,S] r'_i.\]
\end{prop}

\textbf{Proof:}

For $r \in S, v \in \aut(S)$ we have $[v,S]r= r^v[v,S]$. By collecting the elements $r_i$ to the left we may therefore assume that $r_i=r'_i=1$.
Put $M_j= [v_1,S] \cdots [v_j,S]$ and note that $M_1 \subset M_2 \subset \cdots \subset M_m$. As $m>|S|$ there is some $j>1$ such that $M_{j-1} [v_j,S]=M_j=M_{j-1}$. This implies that
$M_{j-1} H=M_{j-1}$ where $H=\la [v_j,S] \ra$. We claim that $H=S$ and therefore $M_{j-1}=S$.
To prove the claim choose $g \in S$ such that $g^{v_j} \not = g$ and for an element $s \in S$ consider $[v_j,s] ([v_j,gs])^{-1} \in H$. Since
\[ [v_j,s] ([v_j,gs])^{-1}= s^{-v_j} (g^{-1}g^{v_j}) s^{v_j} \] and $s \in S$ is arbitrary we see that $(g^{-1}g^{v_j})^S \subset H$ and therefore $H=S$. The Proposition follows. $\square$

We summarise some of the arguments in the proof in the following Proposition which we shall use later. 
\begin{prop}\label{observe} Let $G$ be a finite group with a semisimple subgroup $N=S^{(k)}$. For integers $d$ and $m$ let $g_1, \ldots g_d \in G$ generate a transitive group acting on the set $O$ of simple factors of $G$. Assume that either condition (1) or condition (2) below hold. 

Condition (1): Suppose that $C$ is subset of $\aut (S)$ with the property that $S= \prod_{i=1}^m [\a_i,S]$ for any $\a_1, \ldots \a_m \in C$ and assume that 
\[ |\{   (g_i,S) \ |\ 1 \leq i \leq d, S \in O, \ S^{g_i}=S, \ \exists \a \in C \ (u^{g_i}= u^\a  \ \forall u \in S) \ \}| \geq m .\] 

Condition (2): Suppose that there exists $S \in O$, $s \in \N$, automorphisms $\b_1, \ldots, \b_s \in \aut(S)$ and indices $1 \leq i_1<i_2< \cdots  \ldots i_s \leq d$ such that for all $u \in S$ we have $u^{g_{i_j}} =u^{\b_j}$ for $j=1,2,\ldots, s$ and in addition $\prod_{j=1}^s [\b_j,S]=S$. 
\medskip

Then $N= \prod_{i=1}^d [g_i,N]$.

\end{prop}
\subsection{Proof of Proposition \ref{key}}

The proof of Proposition \ref{key} relies on the following Lemma whose proof occupies the next section.

\begin{lem} \label{derived}
 Given a word $w \in F_d$ there are integers $c(w),f(w) \in \mathbb N$ with the following property:

Suppose $G$ is a finite group with a semisimple normal subgroup $N=S^{(k)}$ with $S$ simple and $|S|>c(w)$. Given any $d$-tuples $\mathbf a_i \in G^{(d)}$ for $i=1,2, \ldots, f(w)$ we can find some other $d$-tuples $\mathbf b_i \in G^{(d)}$ such that 

1. $\mathbf a_i \equiv \mathbf b_i$ mod $N^{(d)}$ and 

2. For $g_i= w(\mathbf b_i)$, the map $\phi : N^{c} \rightarrow N$
defined by
\[ \phi(  x_1 , \ldots, x_f)= [g_1,x_1] [g_2,x_2] \cdots [g_f,x_f]  \]
is surjective.
\end{lem}

Assuming this for the moment we can complete
\textbf{Proof of Proposition \ref{key}:} For tuples $\mathbf r_i \in N^{(d})$ and elements $y_i \in N$ we have
\[  \prod_{i=1}^f w(\mathbf a_i \mathbf r_i)^{y_i}= \prod_{i=1}^f w(\mathbf a_i \mathbf r_i)  \cdot \prod_{i=1}^f  [w(\mathbf a^{t_i}_i \mathbf r^{t_i}_i),y^{t_i}_i]\] where
$t_i:= \prod_{j={i+1}}^f w(\mathbf a_j \mathbf r_j)$.   Define 
$\mathbf a'_i :=\mathbf a_i^{\prod_{j={i+1}}^f w(\mathbf a_j)}$ and observe that the map $\xi : (\mathbf a_1, \ldots, \mathbf a_f)N^{(fd)} \rightarrow (\mathbf a'_1, \ldots, \mathbf a'_f)N^{(fd)}$ defined by 
\[ \xi(\mathbf a_1 \mathbf x_1, \ldots, \mathbf a_f \mathbf x_f) =(\mathbf a_1^{t_1} \mathbf x_1^{t_1}, \ldots, \mathbf a_f^{t_f} \mathbf x_f^{t_f})  \quad \mathbf x_i \in N^{(d)}, \ i=1,2, \ldots f\] is a bijection. We can now apply Lemma \ref{derived} to the $d$ tuples $\mathbf a'_1, \ldots, \mathbf a'_f \in G^{(d)}$ to deduce the existence of some 
$\mathbf r_i \in N^{(d)}$ such that if $\mathbf b_i= \mathbf a_i^{t_i} \mathbf{r_i}^{t_i}$ then the map $\phi: N^{(f)} \rightarrow N$
\[ \phi(y_1, \ldots y_f) = \prod_{i=1}^f [w(\mathbf b_i), y_i^{t_i}] \] is surjective.

We set now $\mathbf z_i = \mathbf a_i^{-1}(\mathbf a_i \mathbf r_i)^{y_i}$ so that 
\[ \psi (\mathbf z_1, \ldots, \mathbf z_f)= \prod_{i=1}^f w(\mathbf a_i \mathbf r_i)  \cdot \prod_{i=1}^f  [w(\mathbf a^{t_i}_i \mathbf r^{t_i}_i),y^{t_i}_i]\left( \prod_{i=1}^f w(\mathbf a_i)  \right)^{-1}=
\]
\[ = \prod_{i=1}^f w(\mathbf a_i \mathbf r_i)  \cdot \phi (y_1, \ldots, y_f) \cdot \left( \prod_{i=1}^f w(\mathbf a_i)  \right)^{-1}.\]
The surjectivity ot $\phi$ now implies the surjectivity of $\psi$.
 $\square$

\subsection{Proof of Lemma \ref{derived}.}

\subsubsection{Preliminaries}
  
Let $r(G)$ denote the minimal degree of a non-trivial real representation of a finite group $G$. We refer to the following lemma as 'the Gowers trick'. 
\begin{lem}[\cite{gow}] \label{gowers} Suppose that $X_1, \ldots, X_n$ are $n \geq 3$ subsets of a finite group $G$ each of size at least $|G|r(G)^{-\frac{n-2}{n}}$. Then $X_1 \cdots X_n=G$. 
\end{lem}

The following proposition summarises some basic results about the sets $[\a, G]$.
\begin{prop} \label{basic} For automorphisms $\a_i$ of a group $G$ and element $g \in G$ we have
\begin{enumerate} \item $[\a_n \a_{n-1} \cdots \a_1, G] \subseteq [\a_1,G] [\a_2, G] \cdots [\a_n,G]$,
\item $[\a,G] g = g^\a [\a, G]$,
\item $[\a,G]^G= \{ h^{-1} [\a,G] h \ | \ h \in G \} \subseteq [\a, G][\a^{-1},G]$ and
\item If $\a$ is the inner automorphism of $G$ given by conjugation by $g$ then $[\a,G]= g^{-1} g^G$. 
\end{enumerate}

\end{prop}

The following Lemma is a special case of \cite{DanN12}, Lemma 4.25. Similar version holds for any universal cover of a simple group of Lie type but we won't need this result in full generality.

\begin{lem} \label{SL} Let $\phi$ be a field automorphism of order $k$ of $S=SL(n,q)$ and let $g \in S$. Put $G= SL(n,q^{1/k})$ and $z=gg^\phi \cdots g^{\phi^{k-1}}$.
Then \[ | \{ a \in S \ | \ a^{g}=a^{\phi} \}| = | C_G(z')| \] where $z' \in G$ is a conjugate to $z^{-1}$ over the algebraic closure $K$ of $\mathbb F_q$.
\end{lem}

\textbf{Proof:} Let $\sigma$ be the Steinberg automorphisms of the algebraic group $SL(n, K)$ induced by the field automorphism $x \mapsto x^q$. We have that $\phi^k=\sigma$. By Lang's theorem there is some $u \in SL(n,K)$ such that $g=uu^{-\phi}$. Put $b=a^u$, then $a^g=a^{\phi}$ is equivalent to $b=b^\phi$ i.e. $b \in G$. Moreover $a=a^\sigma$ implies that 
$b^{u^{-1}}=b^{u^{-\sigma}}$ and so $b$ commutes with $z'=u^{-1}u^\sigma$. Note that $g=uu^{-\phi}=g^\sigma$ implies that $(z')^\phi=z'$ i.e. $z \in G$. Finally observe that $u(z')^{-1}u^{-1}=uu^{-\sigma}=gg^\phi \cdots g^{\phi^{k-1}}=z$.  $\square$

\begin{prop} Let $G$ be a group and $B$ a finite subgroup of $\aut(G)$. Let $n>|B|$ and $\a_1, \ldots, \a_n \in B \backslash \{1\}$. Then
$\prod_{i=1}^n [\a_i,G]$ contains some nontrivial conjugacy class of $G$. 
\end{prop} 
Proof: Observe that for some $1 \leq i<j \leq n$ we have $\a_j \a_{j-1} \cdots \a_i =1$. The result now follows from Proposition \ref{basic} (1) and (3).
\begin{prop} \label{conj} For every integer $l$ there is $d=d(l) \in \N$ such that if $C_1, \ldots, C_d$ are $d$ nontrivial conjugacy classes in a finite simple group of Lie type $S$ of rank at most $l$ then 
$C_1 C_2 \cdots C_d =S$.
\end{prop}
Note that this implies that Proposition \ref{small} holds with $f(l)=d(l)$ whenever all automorphisms $\a_i$ are inner.

The following Proposition can be deduced from Theorem 1 and Lemma 2 from \cite{classicaldecomp}.
\begin{prop} \label{deco} Let $S$ be a classical finite simple group of (untwisted) Lie rank $r>10$. There exists a subgroup $S_0$ isomorphic to $SL(n,q)$ with $n \geq  r/2 - 3$ such that the following hold.

1. $S$ is a product of at most 1000 conjugates of $S_0$,

2. $S_0$ is invariant under $D_S\Phi_S \Gamma_S$,
 
3. For any automorphism $\a \in \aut(S)$ there is an inner automorphism $g$ such that $(\a g)|_{S_0} \in \Phi_{S_0}$.
\end{prop}

\begin{prop} \label{small} There is a function $f_0: \N \rightarrow \N$ with the following property: Let $l \in \N$, let $n=f_0(l)$ and let $S$ be a finite simple group of Lie rank at most $l$. For $i=1, \ldots, n$ let $\a_i \in \aut(S) \backslash \{1\}$. Then $S= [\a_1,S] \cdots [\a_{n},S]$.
\end{prop}

\textbf{Proof:}

Note that $\mathrm{Inn}(S) D\Phi$ is a subgroup of index at most 6 in $\aut(S)$. Using Propostion \ref{basic} (1) implies that each product $[\a_{i},S][\a_{i+1},S] \cdots [\a_{i+5},S]$ contains either a nontrivial conjugacy class or a set $[\b,S]$ where $1 \not = \b \in \mathrm{Inn}(S)D\Phi$. Therefore, putting $f_0(l)=6f_1(l)$ we aim to prove Proposition \ref{small} under the assumption that all automorphisms $\a_i \in \mathrm{Inn}(S) D \Phi$. Let $\a = a \phi$ where $a$ is an inner-diagonal automorphism of $S$ and $\phi \in \Phi$ has order $k$.
The argument in the proof of Lemma \ref{SL} shows that $|C_S(\a)| \leq |C_L(\phi)|$ for the adjoint algebraic group $L$ asociated to $S$. In particular $|C_S(\a)|<|S|^{c/k}$ for some absolute constant $c>0$. Now $|S| < r(S)^C$ for some number $C>0$ depending on $l$ only. We deduce that there exists an integer $N=N(l)$ such that if the order $k$ of $\phi$ is greater than $N$ then $|[\a, S]|>\frac{|S|} {r(S)^{1/3}}$.

Let $M=N!$ and let $\Gamma_0$ be the subgroup of $\Gamma$ of elements of order dividing $M$. Since $\Gamma$ is a cyclic group $|\Gamma_0| \leq M$.
Let $f_1=3+ (l+1)Md(l)$ where $d(l)$ is the integer from Proposition \ref{conj}. If at least three of the automorphisms $\a_1, \ldots, \a_{f_0}$ lie outside $Inn(S)D\Gamma_0$ then we have $[\a_j,S] >\frac{|S|} {r(S)^{1/3}}$ for at least three of the $\a_i$ and we are done Lemma \ref{gowers}. Otherwise at least $(l+1)Md(l)$ of the $\a_1, \ldots, \a_{f_0}$ lie in $Inn(S)D\Gamma_0$. Since $[Inn(S)D \Gamma_0 : Inn(S)] \leq M(l+1)$ repeated applications of Proposition \ref{basic} (1) and (4) give that $\prod_{i=1}^{f_0}[\a_i,S]$ contains a product of some $d$ nontrivial conjugacy classes of $S$ and therefore is all of $S$ by Proposition \ref{conj}. $\square$

Let us denote by $P_n$ the subgroup of permutation matrices in $SL(n,q)$ and by $j_n$ the permutation matrix corresponing to the involution $(1n)(2,n-2)...$ of $S_n$.  Define
\[ A= \left \{ \begin{array}{l} \{ \mathrm{diag} (x,x^{j_m}) \ | \ x \in P_m\} \subset SL(2m,q) \textrm{ if } n=2m \\
\{ \mathrm{diag} ( x,1,x^{j_m}) \ | \ x \in P_m\} \subset SL(2m+1,q) \textrm{ if } n=2m+1 \end{array} \right. \]
Note that $A \simeq \mathrm{Alt}([n/2])$ commutes with all field automorphisms of $SL(n,q)$ and with the graph automorphism $ x \mapsto (x^{-T})^{j_n}$.

\begin{prop} \label{class} For every word $w$ there are integers $l,m$ and $m_0$ with the following property: Let $n>l$ and $S= SL(n,q)$. 

(a) There exists an element $a \in A <S$ such that for any $m_0$ automorphisms $\phi_i \in \Phi \Gamma \subset \aut(S)$ we have \[ S= \prod_{i=1}^{m_0} [a \phi_i,S].\]

(b) There exist word values $g_1, \ldots, g_{m} \in A < S$ such that for any $m$ automorphisms $\phi_i \in \Phi \Gamma \subset \aut(S)$ we have \[ S= \prod_{i=1}^{m} [g_i \phi_i,S].\]
\end{prop}

\textbf{Proof:}
We need the following theorem from \cite{LS}.
\begin{thm} \label{LST} Given a word $w$ there is an integer $k$ such every element of a nonabelian finite simple group of size at least $k$ is a product of two word values.
\end{thm}

The proof of the following Proposition is elementary.
\begin{prop} \label{centraliser}
Let $V$ be a vector space with $V= V_1 \oplus V_2 \oplus \cdots \oplus V_{k+1}$ for subspaces $V_i$ with $\dim V_i=n_i$. Let $A \in \mathrm{End}(V)$ be a matrix which preserves each $V_i$ and acts on $V_i$ as a permutation matrix $a_i$ with $m_i$ cycles. Suppose that $a_1, \ldots, a_k$  have pairwise coprime orders and have no fixed points, while $a_{k+1}=1$. 
Then \[ \dim C_{\mathrm{End}(V)}(A) \leq n_{k+1}^2 + \sum_{i=1}^k m_i(n_i+n_{k+1}). \] 
\end{prop}
Let us first show how (a) implies (b) with $m=2m_0$.
By Theorem \ref{LST} there is some $l=l(w)$ such that if $n>l$ then any element of $A \simeq \mathrm{Alt}([n/2])$ is a product of two word values. Let $g_1,g_2 \in A$ be the word values such that $a=g_2g_1$ and note that by Proposition \ref{basic} $[g_1 \phi_1,S][g_2\phi_2 S] \supseteq [g_2g_1 \phi_2 \phi_1, S]$ since $g_1,g_2$ are fixed by $\phi_1, \phi_2$. 

So it remains to prove part (a). Choose $l>10^3$ so large that there are at least three distinct primes in the interval $(\frac{33n}{200}, \frac{n}{6})$ for each $n>l$. Let $p_1,p_2,p_3$ be such three primes and let $a \in A$ be the permutation matrix with two cycles of lengths $p_1,p_2,p_3,x$ each,  where $x$ is the largest odd integer not exceeding $[n/2] -p_1-p_2-p_3$. Note that $a$ fixes a subspace of $\F _q^n$ of dimension at most $3$. 

We claim now that if the order $v$ of $\phi \in \Phi$ is at least $6$ then $|C_S(a\phi)| < r(S)^{3/4}$. Note that $r(S)= q^{n-1}-1$ when $n>1$ by \cite{LandSeitz}. If $v$ is divisible by at least two of the primes $p_i$ then $v \geq p_ip_j> n^2/7^2> 2n$ and therefore $|C_S(a\phi)| \leq |SL(n,q^{1/v})| < q^{n^2/v} <q^{n/2} \leq r(S)^{3/4}$.

On the other hand if $v$ is divisible by exactly one prime from $p_i$, say by $p_1$  then $a^v$ has two cycles of length $p_2$ and $p_3$ each and has remaining cycles of total length $n-2p_2-2p_3 <(1-33/50)n$. By Proposition \ref{centraliser} the centralizer of $a^v$ in $M_n(F)$ has dimension at most $2(n-2p_3)+ 2(n-2p_2) + (\frac{17}{50}n)^2$ and so 
\[ |C_S(a\phi)|< q^{\frac{4n+ (\frac{17}{50}n)^2}{v}} \leq q^{3(n-1)/4} -1 < r(S)^{3/4}\] since $v>33n/200$ and $n>10^3$.

On the other hand if $v$ is divisible by none of the primes $p_i$ then $a^v$ has at most $2v+3$ cycles of total length at most $n/100$ in addition to the six cycles of length $p_i$. Therefore by Proposition \ref{centraliser} the dimension of $C_{M_n(F)}(a^v)$ is at most $(2v+3)n/100 + \sum_{i=1}^3 2 (2p_i+n/100)$ and again using Lemma \ref{SL} we have

\[  |C_S(a\phi)|< q^{\frac{(2v+3)n +6 \times 35n}{100v}} \leq q^{3(n-1)/4} -1 \leq r(S)^{3/4}\] as $v \geq 6$. 

Notice that if $b=a^{60}\in A$ then $b$ has two cycles of length $p_i$ for each $i=1,2,3$. Therefore since $|b^S|>\sqrt{|S|}$ there is a small constant $c_0$ such that $S$ is a product of $c_0$ conjugacy classes $b^S$.

Let $m_1= 300c_0 + 8$ and suppose that all automorphisms $\phi_i \in \Phi$. I claim that $\prod_{i=1}^{m_1} [a\phi_i,S]=S$. Suppose first that are at least $8$ indices $i$ such that the order of $\phi_i$ is at least 5. We proved above that in this case $|[a\phi_i,S]|> |S| r(S)^{-3/4}$ and the claim follows from Lemma \ref{gowers} with $n=8$. 

Otherwise there will be at least $60c_0$ automorphisms $\phi_i$ which are equal to each other and of order at most 5. Therefore since $[a\phi_i,S]^{60} \supseteq  [a^{60},S]=b^{-1}b^S$ the claim follows from the choice of $c_0$.

Finally we can prove the Proposition in full generality by observing that when $\phi_i \in \Phi \Gamma$ then $[a\phi_1,S][a \phi_2 S]$ contains either $[a\phi,S]$ or $[a^2\phi,S]$ for some $\phi \in \Phi$ and noting that $a^2$ is a conjugate to $a$ in $A$.  $\square$

\subsubsection{Proof of Lemma \ref{derived}}

Let $l=l(w)$ and $m=m(w)$ be the integer in proposition \ref{class}. We will divide the proof into three cases depending on the simple factors of $N$ being 

1. Lie groups of rank at most $10l$,

2. Alternating groups, and 

3. Classical groups of rank at least $10l$. \medskip

\textbf{Case 1: When $S$ is a group of Lie type and rank at most $l$.} Let $c_0=c_0(w)$ be the number from Proposition \ref{coset} and suppose $|S|>c_0$. Take $f>4+4D +4f_0(10l)$ where $f_0=f_0(10l)$ is the function from Proposition \ref{small}. By Proposition \ref{coset} can choose $d$-tuples $\mathbf b_i \in \mathbf a_iN^{(d)}$ such that the word values $g_i:= w(\mathbf b_i)$ do not centralize any simple factor of $S$. We now show that the map $\phi$ is surjective following the same argument from the proof of Theorem \ref{p}.

Without loss of generality we may assume that $\langle g_1, \ldots, g_f \ra$ generates a transitive group on the set $O$ of direct factors of $N$. Let $n_i$ be the number of fixed points of $g_i$ on $O$. If $\sum n_i < f_0$ then the elements $g_i$ satisfy the conditions of Propostion \ref{9.1} and therefore $\phi$ is surjective.
Now assume that $\sum_i n_i \geq f_0$. Proposition \ref{small} gives that $S=\prod_{i=1}^{f_0}[\a_i,S]$ for any nontrivial automorphisms of $S$. We can therefore apply Propostion \ref{observe} with $C= \aut(S) \backslash \{1\}$ in Condition (1) to conclude that $\phi$ is surjective.
\medskip

\textbf{Case 2: When $S$ is an alternating group} 

Let $\mathbf a_i =(a_{i,1}, \ldots, a_{i,d}) \in G^{(d)}$. Let us deal first with the situation when the action of $a_i$ on $N$ is given by permutation of the factors, i.e. when
$(s_1, \ldots ,s_k)^{a_{i,j}}= (s_{\pi_{i,j}(1)}, s_{\pi_{i,j}(2)}, \ldots, s_{\pi_{i,j}(k)})$ for a permutation $\pi_{i,j} \in \mathrm{Sym}(k)$. We note that in this situation whenever $a_{i,j}$ preserves a simple factor of $N$ then it fixes it. Choose now elements $x_1, \ldots, x_d \in S$ such that the conjugacy class $ v^S$ of the word value $v:=w(x_1, \ldots, x_d) \in S$ is such that $(v^S)^3=S$. Define now $\mathbf b_i= \mathbf a_i \mathbf x$ where $\mathbf x= (\bar x_1, \bar x_2, \ldots, \bar x_d) \in N^{(d)}$ with $\bar x_i=(x_i,x_i, \cdots, x_i) \in S^{k}=N$. Note that $S_j \in O$ is any simple factor of $N$ preserved by $w(\mathbf b_i)$ then $g^{w(\mathbf b_i)}=g^v$ for any $g \in S_j$.

It is enough to take $f'= 4+4D +12$ and argue as above noting that if the sum $\sum_{i=1}^f n_i$ of fixed points of $w(\mathbf b_i)$ on $O$ is at least 3 then the surjectivity of $\phi$ follows from $[v,S]^3=v^{-3} (v^S)^3 =S$ and Proposition \ref{observe} with $C= \{ u \mapsto u^v  \ \forall u \in S\}$. 

Now we can deal with the general case without the retrictions on the $a_{i,j}$.  Suppose $S=\mathrm{Alt}\{1,2,\cdots,n\}$ and let $A_1= A_3=A_5=\mathrm{Alt}\{1,2 \ldots,n-2\}$, 
$A_2= A_4=A_6=\mathrm{Alt}\{3,4, \ldots, n\}$. It is easy to see that $S=(A_1A_2)^3=A_1 \cdots A_6$. Define further $N_i= A_i^{(k)} \leq N$.

Let $T_1,T_2$ be the subgroups of order $2$ in $\aut(S)$ generated by conjugation with the transposition $(n-1,n)$ and $(12)$ respectively. Note that $T_1 \mathrm{Inn}(S)= \aut(S)= T_2 \mathrm{Inn}(S)$. Let us put $f=6f'$ and relabel the tuples $\mathbf{a}_i$ as $\mathbf a_{i,j}$ with $i=1, \ldots, 6, j=1, \ldots, f'$. We can of course replace each $\mathbf a_{i,j}$ by any member of its coset $\mathbf a_{i,j} N^{(d)}$ and in this way we can ensure that for $i=1,3,5$ the constutuents of each $\mathbf a_{i,j}$ act on $N$ as elements of $T_1 \wr \mathrm{Sym}(k)$ while for $i=2,4,6$ they acts as elements of $T_2 \wr \mathrm{Sym}(k)$.

This ensures that each $\mathbf a_{i,j}$ acts on $N_i$ by permutation of its simple factors $A_i$.
By the special case above we can find $d$-tuples $\mathbf b_{i,j} \equiv \mathbf a_{i,j}$ mod $N$ such that the maps $\phi_i : N_i^{(f')} \rightarrow N_i$
 \[ \phi_{i}(y_1, \ldots, y_{f'})= \prod_{j=1}^{f'} [w(\mathbf a_{i,j}), y_j] \quad y_i \in N_i\] are surjective. Therefore $\phi= \phi_1 \cdots \phi_6$ is surjective onto $N$ and this case is completed.

\textbf{Case 3: When $S$ is a classical group of rank at least $10 l(w)$.}

Proposition \ref{deco} gives that $S=A_1 \cdots A_{200}$ where each $A_i$ is isomorphic to $SL(n,q)$ with $n>l(w)$ and invariant under $D \Phi \Gamma$.
Following the same argument as in the previous case we can set $f>10^3 f''$ provided we prove the Lemma with $f''$ in the case $S=SL(n,q)$ and the automorphisms $\mathbf a_i$ act on $N=S^{(k)}$ as elements of $\Phi \Gamma \wr \mathrm{Sym}(k)$. 

As before let $n_i$ be the number of fixed points of $w(\mathbf b_i)$ on the set $O$ of simple factors of $N$. Define $c_i$ to be the number of orbits of $w(\mathbf a_i)$ on $O$. Since $w(\mathbf a_i)$ and $w(\mathbf a_i)$ act in the same way on $O$ we have $c_i \leq (k+n_i)/2$.

Suppose first that $\sum_{i=1}^{f''} n_i \leq (1-m^{-1})kf''$. Then \[ \sum_{i=1}^{f''} c_i \leq \sum_{i=1}^{f''} (k+n_i)/2 \leq kf'' - \frac{kf''}{2m} < kf'' -2k-2D,\]
since $f''>m(4+4D)$ and we may apply Proposition \ref{9.1} to deduce that $\phi$ is surjective for any choice of $\mathbf b_i \equiv \mathbf a_i$ mod $N$.

Now suppose that $\sum_{i=1}^{f''} n_i \leq (1-m^{-1})kf''$. This implies the existence of some simple factor $S \in O$ such that the set of indices 
\[ P :=\{ 1 \leq i \leq f'' \ | \  S^{w(\mathbf a_i)}=S \}\] has cardinality  $|P| \geq f''(1-m^{-1})$. Without loss of generality we may assume that $S=S_1$. Now divide the interval $[1,f'']$ into $m$ equal consecutive subintervals $I_1, \ldots I_m$. As $|P| \geq f''(1-m^{-1})$ we deduce that $P \cap I_s \not = \emptyset$ for every $s \in \{1, \ldots, m\}$.

 Let $g_1,\ldots, g_m$ be the word values in the subgroup $A$ provided by Proposition \ref{class} (2). For $s =1, \ldots, m$ suppose that $g_s=w(v_{s,1}, \ldots, v_{s,d})$. Define
$\mathbf v_s = ( \bar v_{s,1}, \bar v_{s,2}, \cdots , \bar v_{s,d}) \in N^{(d)}$ where $\bar v_{s,i}= (v_{s,i}, v_{s,i}, \cdots ,v_{s,i}) \in S^{(k)}=N$. Note that for any simple factor $S \in O$ we have $u^{w(\mathbf v_s)} =u^{g_s}$ for all $u \in S$.

For any index $i \in I_s$ we define $\mathbf x_i:= \mathbf v_s$, $\mathbf b_i= \mathbf a_i \mathbf x_i$ and observe that if $i \in P$ and $u \in S=S_1$ then $u^{w(\mathbf b_i)}=u^{g_s \phi_i}$ where $\phi_i \in \Phi \Gamma$ is the action of $\mathbf a_i$ on $S$.

Now Proposition \ref{class} ensures that $S= \prod_{s=1}^m [g_i \phi_i,S]$.  The sequence of sets $\{ [w(\mathbf b_i),S] \ | \ i \in P \}$ contains at least one occurence of $[g_s \phi_s,S]$ for every $1 \leq s \leq m$ because $P \cap I_s \not = \emptyset$. Therefore the map $\phi$ is surjective by Proposition \ref{observe}. The proof of Lemma \ref{derived} is complete.

\subsection{Proof of Proposition \ref{coset}}

We will need several standard results about polynomial rings.

 For a commutative ring $R$ let $u_+(a)$ ( respectively $u_-(a)$) denote the upper ( respectively lower) unitriangular 2 by 2 matrix with entry $a \in R$ in $SL(2,R)$. The following Proposition is a standard application of the Ping-Pong Lemma, cf \cite{delaHarpe}[II.B, Lemma 24].

\begin{prop} Let $R$ be a polynomial ring and let $I$ be the ideal of $R$ consisting of polynomials with constant term 0.  Let $U_+(I) = \{ u_+(a) \ | \ a \in I \} \leq SL(2,R)$ and $U_-(I) = \{ u_-(a) \ | \ a \in I \}$. The group generated by $U_+(I)$ and $U_-(I)$ in $PSL(2,R)$ is isomorphic to the free product $U_+(I) \star U_-(I)$. 
\end{prop}

For a field $F$ and a set $X$ we denote by $F[X]$ the polynomial ring in commuting variables from $X$. We will say that two non-zeoro polynomials $f,g \in F[X]$ are \emph{disjoint} if they have disjoint variables, i. e, if $f \in F[X_1]$ and $g \in F[X_2]$ for disjoint subsets $X_1,X_2 \subset X$. More generally two finite subsets $P_1,P_2 \subset R \backslash \{0\}$ of polynomials are called disjoint if 
there are disjoint subsets $X_1,X_2 \subset X$ such that $P_i \subset F[X_i]$ for $i=1,2$. This is equivalent to the condition that $\prod_{f \in P_1}$ is disjoint from $\prod_{f \in P_2}$.  
\begin{cor} \label{free} Suppose $R=F[X]$ with ideal $I$ as above and for any nonzero polynomials $f_{1,s},f_{2,s},f_{3,s} \in I \vartriangleleft F[X]$ with  $s=1,\ldots, k$ let \[ x_i= u_-(f_{1,i})u_+(f_{2,i})u_-(f_{3,i}) \in SL(2,R).\] Suppose that $f_{3,s}$ and $f_{1,s+1}$ are disjoint for all $s=1, \ldots, k-1$. Then $x_1x_2 \cdots x_k$ is not central in $SL(2,R)$.
\end{cor}

\begin{lem} \label{field} Let $f \in \F [t_1, \ldots, t_n]$ be a non-zero polynomial with coefficients in a finite field $\F$, such that the degree of $f$ in each variable $t_i$ does not exceed $N \in \N$. Assuming that $\F_0$ is a subfield of $\F$ such that $|\F_0|>N$ there exist $y_1, \ldots, y_n \in \F_0$ such that $f(y_1, \ldots, y_n)  \not = 0$.
\end{lem}

\textbf{Proof of Proposition \ref{coset}.}

We can replace the $d$-tuple $\mathbf g =(g_i)$ with any $\mathbf g \mathbf n$ where $\mathbf n  \in N^{(d)}$. Suppose first that $S= \mathrm{Alt}(n)$ is an alternating group. We may assume that each $g_i$ acts on $N=S^{(m)}$ as an element of $C \wr \mathrm{Sym}(m)$ where $C < \aut(S)$ is the group generated by conjugation with fixed transposition $t \in S$. Let $A_0$ be the subgroup of $C_S(t)$ ismorphic to $\mathrm{Alt}(n-2)$. For an element $v \in A_0 <S$ let us write $\bar v = (v,v, \ldots, v) \in S^{(m)}=N$. It is clear that $\bar v \in A_0^{(m)}$ is centralized by each $g_i$ and therefore for any $v_i \in A$ we have $w(g_1 \bar v_1, \cdots, g_d \bar v_d)= w(g_1, \ldots g_d) \bar z$ where $z=w(v_1, \ldots, v_d) \in A_0$. Now if $n$ is sufficiently large compared to $w$ we can find $v_i \in A_0$ such that $z \not =1$. Finally note that if $w(g_1 \bar v_1, \cdots, g_d \bar v_d)$ preserves some factor $S$ of $N$ then it acts on it as conjugation by either $z$ or $zt$ and so cannot centralize this factor.

 We will now deal with the main case when $S$ is of Lie type.  In this case we may assume that each $g_i$ acts on $N= S^{(m)}$ as an element of the subgroup $T=(D \Phi \Gamma) \wr \mathrm{Sym}(m) \subset \aut(S) \wr \mathrm{Sym}(m) =\aut(N)$. If $S$ has large Lie rank as some function of $w$ then Proposition \ref{deco} gives a subgroup $S_0=PSL(n,F)$ invariant under $D\Phi \Gamma$ and therefore $S_0^{(m)} \leq N$ is invariant under $T$. By restricing the elements $a_i$ to vary over $S_0^{(m)}$ we may assume that $S=S_0$ and the rank $n-1$ of $S_0$ is large in terms of $w$. Now use the argument above with the subgroup $A$ from Proposition \ref{class} in place of $A_0$.

It therefore remains to deal with the case when the rank of $S$ is bounded in terms of $w$. Therefore the field of definition $\F$ of $S$ can be assumed to be sufficiently large.

We can make a further reduction. If the type of $S$ is not $G_2,B_2, ^2 B_2$ or $F_4$ then $S$ has a quasisimple subgroup $S_1$ of type $A_1$ which is preserved by $D \Phi \Gamma$.
When $S$ has type $G_2$ or $F_4$ then $S$ has  $D\Phi \Gamma$- invariant semisimple subgroup $S_1$ of type $A_1$ if $\Gamma = \{1\}$, or of type $A_1 \times A_1$ if $S$ has exceptional graph automorphisms (which exist only in characteristic 3 for $G_2$ and characteristic 2 for $F_4$).  When $S=B_2(q)$ and $q$ is odd then $\Gamma =1$ and therefore there is a $D\Phi$-invariant subgroup $S_1$ of type $A_1$. Finally if $q$ is even then every diagonal automorphism of $S$ is inner and we can take $S_1= \ ^2B_2(q) < S$ which is invariant under $\Phi \Gamma$. By restricting $a_i \in N$ to range over $S_1^{(m)}$ we may now reduce to the case when $S=S_1$ is of type $A_1$ or $^2B_2$.
Moreover we may continue to assume that $g_i \in T_0:=(D \Phi)_S \wr \mathrm{Sym}(m) < \aut(N)$. \medskip

\textbf{I: When $S=PSL(2,\F)$.} 

We need to introduce some notation.
Suppose $w = \prod_{s=1}^k x_{i_s}^{\e_{s}}$ where $e_{s} \in \{\pm 1 \}$ and $i_1, \ldots, i_k \in \{1, \ldots, d\}$. We have 
\[ w(g_1a_1, \ldots, g_da_d)= w(g_1, \ldots, g_d) \prod_{s=1}^k a_{i_s}^{\e_{s}f_{s}} \] for some automorphisms $f_{s} \in T_0 < \aut(N)$. Set $w_0:=w(g_1, \ldots g_d)$ and let $J \subseteq \{1, \ldots, m\}$ be the set of indices $j$ such that $S_j^{w_0}=S_j$ and such that $w_0$ induces an the inner automorphism on each $S_j$. For $j \in J$ let $h_j \in S_j$ be the element such that $u^{w_0}=u^{h_j^{-1}} \ \forall u \in S_j$. Note that $w(g_1a_1, \ldots,g_da_d)$ cannot centralize any $S_j$ with $j \not \in J$ for any choice of $a_i \in N$.

Let us write $a_i= (a_{i,1}, \ldots a_{i,m}) \in N$ with $a_{i,j} \in S_j$ and let $f_{s}$ act on $N$ by \[ (\tau_1, \ldots, \tau_m)^{f_{s}}=(\tau_{\pi_s(1)}^{r_{s,1}},\tau_{\pi_s(2)}^{r_{s,2}}, \ldots, \tau_{\pi_s(m)}^{r_{s,m}}), \quad \tau_i \in S, i=1,2,\ldots,m \] for some permutations $\pi_s \in \mathrm{Sym}(m)$ and automorphisms $r_{s,j} \in D \Phi < \aut(S)$.

Denote by \begin{equation} \label{u} u_j= \prod_{s=1}^k a_{i_s, \pi_s(j)}^{\epsilon_s r_{s,j}} \in S_j \end{equation}
 the projection of the element $\prod_{s=1}^k a_{i_s}^{\e_{s}f_{s}}$ onto $S_j$. We need to prove the existence of some $a_1, \ldots a_d \in N$ such that for every $j \in J$ we have $h_j \not = u_j$. We will prefer to work with $a_{i,j}$, $u_j$ and $h_j$ taking values in the universal covering group $SL(2,\F)$ of $S$ and treat this problem as a system of polynomial inequations in the matrix coefficients of $a_{s,j} \in SL(2,q)$. Note that $h_j \not =u_j$ in $S$ becomes equivalent to proving that $h_j^{-1}u_j$ is not central in $SL(2,\F)$. 

\textbf{Subcase I(a):} Suppose that the characteristic $p$ of $\F$ is large compared to the length $k$ of $w$. More precisely let us assume that $p>3k$. \medskip

For $\a=1 \ldots, d$, $\b=1, \ldots, m$ and $l=1,2,3$ let $R$ be the polynomial ring $\F[\{t_{\a,\b}\}]$ in $3dm$ commuting indeterminates $t_{\a,\b,l}$. 
Put  \[ A_{\a,\b}=A_{\a,\b}(\mathbf t_{\a,\b}):= u_-(t_{\a,\b,1})u_+(t_{\a,\b,2})u_-(t_{\a,\b,3})  \in SL(2,R),  \] where
$\mathbf t_{\a,\b}=( t_{\a,\b,1},t_{\a,\b,2},t_{\a,\b,3})$.  In this subcase we are going to search for values $y_{\a,\b,l} \in \F_p$ and specialise the indeterminates $t_{\a,\b,l}\mapsto y_{\a,\b,i}$ to lie in $\F_p$, so that each $a_{\a,\b}:=A( \mathbf y_{\a,\b})$ is going to be fixed under any field automorphism of $\F$. Under this restriction we define each $r \in D\Phi \subset \aut(S)$ to acts on $A_{\a,\b}$ as a diagonal automorphism, in particular there are elements $\lambda^{\a,\b,r}_i \in \F$ which do not depend on $\mathbf t_{\a,\b}$ such that \[ (A_{\a,\b})^{r}= u_-(\lambda_{1}^{\a,\b,r}t_{\a,\b,1})u_+(\lambda_2^{\a,\b,r}t_{\a,\b,2})u_-(\lambda_3^{\a,\b,r}t_{\a,\b,3}).\] We observe that the coefficients of $A_{\a,\b}^{\pm r}$ are polynomials of degree at most 1 in each of $t_{\a,\b,l}$. Put
\[ U_j= \prod_{s=1}^k A_{i_s, \pi_s(j)}^{\epsilon_s r_{s,j}} \in SL(2,R).\]
For each $ \alpha \in \{1, \ldots, d\}$ and $\beta \in \{1, \ldots, m\}$ define
\begin{equation} \label{C} C_{\a,\b}=\{ (j,s) \ | \ 1 \leq j \leq m, \ 1 \leq s \leq k , \ (\alpha,\beta)= (i_s,\pi_s(j)) \}.\end{equation}

The set $C_{\a,\b}$ simply lists all all pairs $(s,j)$ such that $A_{i_s, \pi_s(j)}$ from (\ref{u}) is a matrix whose coefficients are polynomials in $\F[\mathbf t_{\a,\b}]$. Note that each coefficient has degree at most 1 in each of $t_{\a,\b,l}$.
\begin{prop} \label{degree} For each $ \alpha \in \{1, \ldots, d\}$ and $\beta \in \{1, \ldots, m\}$ we have $|C_{\a,\b}| \leq k$.
\end{prop} 

\textbf{Proof:} Given $\alpha, \beta$ there are at most $k$ choices for $s$ such that $i_s= \alpha$ and then $j$ is determined from $\pi_s(j)=\beta$.
$\square$

Now consider the element $h_j^{-1} U_j \in SL(2,R)$. We claim that $h_j^{-1} U_j$ is not a central element of $SL(2,R)$. 

Suppose for the sake of contradiction that $h_j^{-1}U_j$ is central in $SL(2,R)$. We can set all indeterminates $t_{\a,\b,l}$ of $R$ to be $0$,  when clearly $U_j$ vealuates to $1 \in SL(2,\F)$ which implies that $h_j$ is in the centre of $SL(2,\F)$. So actually we must have that $U_j$ is a central element of $SL(2,R)$. 
However it is immediate to check that the matrices $x_s=A_{i_s, \pi_s(j)}^{\epsilon_s r_{s,j}}$ satisfy the conditions in Corollary \ref{free}  therefore $U_j$ cannot be central. This proves the claim.

Suppose that $\left ( \begin{array}{cc} a & b \\ c & d \end{array} \right )$ is the matrix
representing $h_j^{-1}U_j$ in $SL(2,R)$. Since $h_j^{-1}U_j$ is not a scalar, at least one of the polynomials $b,c,a-d$ is not zero, let us denote this non-zero polynomial by $f_j \in R$.  Note that since $h_j \in SL(2,\F)$ each $f_j$ is a  $\F$-linear combination of the coefficients of $U_j=\prod_{s=1}^k A_{i_s, \pi_j(s)}^{\epsilon_s r_{s,j}}$
and moreover each $A_{i_s, \pi_j(s)}^{\epsilon_s r_{s,j}}$ has matrix coefficients from $\F[\mathbf t_{i_s,\pi_s(j)}]$ having degree at most 1 in each of the three variables of $t_{i_s,\pi_s(j),l}$, $l=1,2,3$. Consider now the non-zero  polynomial $f=\prod_{j\in J} f_j \in R$. The above observation and Proposition \ref{degree} show that the maximal degree of any variable $t_{\a,\b,l}$ ocurring in $f$ is at most $k$. Therefore since $p>k$ Proposition \ref{field} provides elements $y_{\a,\b,l} \in \F_p$ such that $f(\mathbf y) \not =0$. Set $a_{\a,\b}:= A_{\a,\b}(\mathbf y_{\a,\b}) \in SL(2,\F)$ and let $u_j$ be defined from (\ref{u}).

Out choice of $y_{\a,\b,l}$ ensures that $u_j \not = h_j$ in $PSL(2,\F)$ for every $j \in J$ and therefore $w(g_1a_1, \ldots, g_d,a_d)$ does not centralize the factors $S_j$ of $N$.
\medskip

\textbf{Subcase I(b):} We now examine the case when $|F|=p^L$ with $L$ large, more precisely assume that $L>k^2$. For $1 \leq s \leq k$ and $1 \leq j \leq m$ let $n_{s,j}$ be defined as the integer from $\{0,1,\ldots, L-1\}$ such that $r_{s,j} \in D \Phi$ has field component $[p^{n_{s,j}}]$. We extend the automorphism $[p^n]$ of $\F$ to an endomorphism of $R$ with the same name defined by $f \mapsto f^{p^n}$ for each $f \in R$. 

Let us assume first assume that \[  \quad n_{s,j} \leq L-k \quad \textrm{ for all } \ s=1, \ldots, k, \ j=1,\ldots, m. \] 

Then arguing exactly as in Case I(a) we deduce that every $A_{i_s, \pi_j(s)}^{\epsilon_s r_{s,j}}$ is a matrix with coefficients which are polynomials from $\F[ \mathbf t_{i_s,\pi_s(j)}]$ of degree at most $n_{s,j} \leq p^{L-2k}$ in each of the variables $t_{i_s,\pi_s(j),l}$. \medskip

 As before since $|C_{\a,\b}| \leq k$ we deduce that the non-zero polynomial $f \in R$ defined in the same way as above has degree at most $kp^{L-k}$ in any variable $t_{\a,\b,l}$. Proposition \ref {field} and the fact that $kp^{L-k} \leq k2^{-k}p^L < p^L=|\F|$ for any $k \in \N$ give the existence of elements $y_{\a,\b,l} \in F$ such that $f(\mathbf y) \not =0$ and again we set $a_{\a,\b}:=A_{\a,\b}(\mathbf y_{\a,\b})$. \medskip

To deal with the general case define $V \subset \aut(\F)$ by \[ V = \{ [p^n]: \F \rightarrow \F \ | \ L-k<n<L\}.\] Given any $\a,\b$ recall that $C_{\a,\b}$ defined in (\ref{C}) is the number of pairs $(s,j)$ such that $a_{i_s, \pi_j(s)}^{\epsilon_s r_{s,j}} \in SL(2,R)$ has coefficients which are polynomials in $\F[\mathbf t_{\a,\b}^{n_{s,j}}]$, moreover these polynomials have degree at most 1 in each of $t_{\a,\b,l}^{n_{s,j}}$. Given any $(s_0,j_0) \in C_{\a,\b}$ the number of integers $0 \leq n<L$ such that $[p^{n+n_{s_0,j_0}}] \in V$ is at most $k$. Since $|C_{\a,\b}| \leq k$ and $L> k^2$ by the pigeonhole principle we concude that there is some integer, denoted $n'_{\a,\b} \in \{0,1, \ldots, L-1\}$ such that $[p^{n'_{\a,\b}+n_{s,j}}] \not \in V$ for any $(s,j) \in C_{\a,\b}$. Now define 
 \[ A'_{\a,\b}:=A_{\a,\b}(\mathbf t_{\a,\b}^{n'_{\a,\b}})= u_-(t_{\a,\b,1}^{p^{n'_{\a,\b}}})u_+(t_{\a,\b,2}^{p^{n'_{\a,\b}}})u_-(t_{\a,\b,3}^{p^{n'_{\a,\b}}})  \in SL(2,R). \] 

For each $(s,j) \in C_{\a,\b}$ we have 
$(A'_{\a,\b})^{\epsilon_s r_{s,j}}$ is a matrix with coefficients which are polynomials from $\F [\mathbf t_{\a,\b}^{p^{n'_{\a,\b}+{n_{s,j}}}}]$ of degree at most 1 in each of  $t_{\a,\b,l}^{p^{n'_{\a,\b}+{n_{s,j}}}}$. From the choice of $n'_{\a,\b}$ it folows that there is an integer $n''_{s,j} \in \{ 1, \ldots, L-2k\}$ such that $x^{p^{n''_{s,j}}}= x^{p^{n'_{\a,\b}+{n_{s,j}}}}$ for every $x \in \F$. Now replace each occurrence of $t_{\a,\b,l}^{p^{n'_{\a,\b}+{n_{s,j}}}}$ in the coefficients of $(A'_{\a,\b})^{\epsilon_s r_{s,j}}$ by $t_{\a,\b,l}^{p^{n''_{s,j}}}$ obtaining the matrix 
$\tilde A_{s,j}(\mathbf t_{\a,\b})$.  Each coefficient of $\tilde A_{s,j}$ is now a polynomial in $\mathbf t_{\a,\b}$ of degree at most $p^{L-k}$ in each variable $t_{\a,\b,l}$. 
For every $\mathbf y_0=(y_1,y_2,y_3) \in \F^{(3)}$ we have $\tilde A_{s,j}(\mathbf y_0)= (A'_{\a,\b}(\mathbf y_0))^{\epsilon_s r_{s,j}}=(A'_{i_s, \pi_j(s)}(\mathbf y_0))^{\epsilon_s r_{s,j}}$. Moreover just as in Case I(a) Proposition \ref{free} gives that $\prod_{s=1}^k \tilde A_{s,j} \in SL(2,R)$ is noncentral. Again we find $y_{\a,\b,l} \in \F$ such that 
\[ h_j^{-1} \prod_{s=1}^k \tilde A_{s,j}(\mathbf y_{i_s, \pi_j(s)}) =h_j^{-1} U_j(\mathbf y), \quad \mathbf y= ( y_{\a,\b,l})_{\a,\b,l} \in \F^{3dm}\] is not in the centre of $SL(2,\F)$ for every $j \in J$ and we set $a_{\a,\b}= A'_{\a,\b}(\mathbf y_{\a,\b})$.
\medskip

\textbf{II: When $S$ is the Suzuki group $^2B_2(2^{2L-1})$.}

We will need a version of Proposition \ref{free} which will apply to a specific linear $4$-dimensional linear representation of $^2B_2$.

Let $R$ be a polynomial ring and as usual let $\deg f$ be the maximal degree of the monomials of $f \in R$. Let $I$ be the ideal of polynomials with constant term zero. Define the subset $D_+ \subset SL(4,R)$ to consist of those upper unitriangular matrices  $g= (g_{i,j})$ with $\deg g_{1,4} > \max \{ \deg g_{i,j} \ | \ (i,j) \not = (1,4) \}$, and set $D_-=D_+^T$.  
\begin{prop} \label{suzfree} Let $U_1$ and $U_2$ are two subgroups of $SL(4,R)$ such that $U_1 \subset D_+ \cup \{1\}$ and $U_2 \subset D_- \cup \{1\}$. Then $U_1$ and $U_2$ generate their free product $U_1 \star U_2$ in $PSL(4,R)$.
\end{prop}

\textbf{Proof:} Consider the natural action of $PSL(4,R)$ on the projective space $\mathbb P(R^4)$ and define two subsets $V_+,V_- \subset \mathbb P(R^4)$ as follows:
\[ V_+= \{ [f_1,f_2,f_3,f_4] \in \mathbb P(R^4) \ | \ \deg f_4 > \max \{ \deg f_1, \deg f_2 , \deg f_3 \}\}, \]
\[ V_-= \{ [f_1,f_2,f_3,f_4] \in \mathbb P(R^4) \ | \ \deg f_1 > \max \{ \deg f_2, \deg f_2 , \deg f_4 \}\}. \]
It is immediate that $g (V_+) \subset V_-$ for all $g \in D_+$ and  $g (V_-) \subset V_+$ for all $g \in D_-$. The proposition follows from the Ping-Pong Lemma \cite{delaHarpe}[II.B, Lemma 24].
$\square$

Suppose now that $S= \ ^2B_2(2^{2L-1})$. Set $\F=\F_{2^{2L-1}}$ define $\theta= 2^L$ noting that $[\theta]^2=[2] \in \aut(\F)$. We assume that $L>4k^2$.  
Define
\[ e_-(t_1,t_2, v_1,v_2)= \left ( \begin{array}{cccc} 1 & 0 & 0 & 0 \\  t_1 & 1 & 0 & 0 \\ v_1+ t_1t_2 & t_2 & 1 & 0 \\ t_1v_1+v_2 + t_1^2 t_2 & v_1 & t_1 & 1 \end{array}\right)\] and set $e_+(t_1,t_2,v_1,v_2)= u_-(t_1,t_2,v_1,v_2)^T$.
We have that $\{ e_-(\tau,\tau^\theta,\nu, \nu^\theta) \ | \ \tau, \nu \in \F\}$ is a parametrisation of a Sylow $2$-subgroup of $S$.

\begin{lem} \label{uv} Let $R$ be a polynomial ring and let $f_1,f_2,g_1, g_2 \in I$ be non-zero polynomials such that $\{f_1,f_2\}$ is disjoint from $\{g_1,g_2\}$.  Then \[ \deg (f_1g_1+g_2 + f_1^2 f_2) = \max \{ \deg (f_1g_1), \ \deg g_2, \ \deg (f_1^2f_2)\} >  \]
\[ > \max \{ \deg f_1, \deg f_2, \deg g_1,\deg (f_1f_2) \}.\]

As a consequence we have $e_-(f_1,f_2,g_1, g_2)$ belongs to $D_- \subset SL(4,R)$ and $e_-(f_1,f_2,g_1,g_2)$ belongs to $D_- \subset SL(4,R)$.
\end{lem}

Lemma \ref{uv} and Proposition \ref{suzfree} together give the following Corollary.

\begin{cor} \label{suzapplic} Let $R$ be a polynomial ring and for $i=1, \ldots, k$, $j=1,2$ and $l=1,\ldots, 3$ let $f_{i,j.l},g_{i,j,l} \in I$ be nonzero polynomials with the following properties:

1. For every $1 \leq i \leq k$ and  $1 \leq l \leq 3$  the set $\{f_{i,1,l},f_{i,2,l}\}$ is disjoint from $\{g_{i,1,l},g_{i,2,l}\}$.

2. For every $i=1, \ldots, k-1$ the set $\cup_{j=1}^2\{f_{i+1,j,1},g_{i+1,j,1}\}$ is disjoint from $\cup_{j=1}^2\{f_{i,j,3},g_{i,j,3}\}$.

Define \[ x_k= e_-(\mathbf f_{i,1},\mathbf{g}_{i,1})e_-(\mathbf f_{i,2},\mathbf{g}_{i,2})e_-(\mathbf f_{i,3},\mathbf{g}_{i,3}) \]
with $\mathbf f_{i,l}= (f_{i,1,l},f_{i,2,l})$ and $\mathbf g_{i,l}= (g_{i,1,l},g_{i,2,l})$.

Then $x_1\cdots x_k$ is not central in $SL(4,R)$.
\end{cor}

Indeed condition (2) ensures that in the product $x_1 \cdots x_k$ there is no cancellation between the last term $e_-(\mathbf f_{i,3},\mathbf{g}_{i,3})$ of $x_i$ and the first term $e_-(\mathbf f_{i+1,1},\mathbf{g}_{i+1,1})$ of $x_{i+1}$.
\medskip

For $\a =1, \ldots, d, \b =1, \ldots,m$ and $l=1,2,3$ let $R$ be the polynomial ring over $\F$ in $6dm$ distinct variables $t_{\a,\b,l},v_{\a,\b,l}$. As before let $n_{s,j}$ be the integers such that $[2^{n_{s,j}}]$ is the field automorphism component of $r_{s,j} \in D\Phi \subset \aut(S)$ in the expression (\ref{u}).

Define $V \subset \aut(\F)$ by \[ V = \{ [p^n]: \F \rightarrow \F \ | \ 2L-2k<n<2L-1\}\] and for $1 \leq \a \leq d, 1 \leq \b \leq m$ let $C_{\a,\b}$ be the set defined in (\ref{C}). 

Given any $(s,j) \in C_{\a,\b}$ the number of integers $0 \leq n<2L-1$ such that $[2^{n+n_{s,j}}] \in V$ or $[2^{n+n_{s,j}+L}] \in V$ is at most $2k+2k=4k$. Since $|C_{\a,\b}| \leq k$ and $L> 4k^2$ by the pigeonhole principle we concude that there is some integer, denoted $n'_{\a,\b} \in \{0,1, \ldots, L-1\}$ such that $[2^{n'_{\a,\b}+n_{s,j}}] \not \in V$ and $[2^{n'_{\a,\b}+n_{s,j}+L}] \not \in V$ for any $(s,j) \in C_{\a,\b}$. \medskip
 
Let \[ E_{\pm}(t,v):= e_{\pm}(t,t^\theta,v,v^\theta) \in SL(4,\F[t,v])\] and define 
\[E_{\pm, \a,\b,l} = E_{\pm}(t_{\a,\b,l}^{2^{n'_{\a,\b}}},v_{\a,\b,l}^{2^{n'_{\a,\b}}}) \in SL(4,R)\] 
\[ A_{\a,\b}=A_{\a,\b}(\mathbf t_{\a,\b}, \mathbf v_{\a,\b}):= E_{-,\a,\b,1}E_{+,\a,\b,2}E_{-,\a,\b,3} \in SL(4,R),\]
where $\mathbf t_{\a,\b}=(t_{\a,\b,1},t_{\a,\b,2},t_{\a,\b,3})$ and $\mathbf v_{\a,\b}=(v_{\a,\b,1},v_{\a,\b,2},v_{\a,\b,3})$. 
Observe that each coefficient of $E_{\pm, \a,\b,l}$ is can be written as polynomial over $\F$ of degree at most two in each of variables 
\[ t_{\a,\b,l}^{2^{n'_{\a,\b}}},\ t_{\a,\b,l}^{2^{n'_{\a,\b}+L}}, \ v_{\a,\b,l}^{2^{n'_{\a,\b}}}, \ v_{\a,\b,l}^{2^{n'_{\a,\b}+L}}. \]
\medskip

As before we consider for $j \in J$ \[ U_j=\prod_{s=1}^k A_{i_s, \pi_s(j)}^{\epsilon_s r_{s,j}} \in SL(4,R)\] and we are looking for field elements $y_{\a,\b,l},z_{\a,\b,l} \in \F$ such that if $u_j$ is the evaluation of $U_j$ at $t_{\a,\b,l} \mapsto y_{\a,\b,l}$, 
$v_{\a,\b,l} \mapsto z_{\a,\b,l}$ then $h_j \not = u_j$. Since $Z(Sp(4,\F))=1$ it is enough to  prove that $h_j^{-1}u_j$ is not central as an element of $SL(4,\F)$.

It is immediate that each $E_{\pm, i_s, \pi_s(j)}^{\pm r_{s,j}}$ is a matrix in $SL(4,R)$ with coefficients which are polynomials of degree at most two in terms of the four monomials 
\[ K_{\a,\b,l}:= \left \{ t_{\a,\b,l}^{2^{n'_{\a,\b}+n_{s,j}+b}}, \ v_{\a,\b,l}^{2^{n'_{\a,\b}+n_{s,j}+b}}\ | \  \ b\in \{0,L\} \right \}. \]

 From the choice of $V$ and $n'_{\a,\b}$ it follows that there are integers $m_{\a,\b},k_{\a,\b}$ from $\{0,1, \ldots 2L-1-2k\}$ such that $x^{2^{n'_{\a,\b}+n_{s,j}}}=x^{2^{m_{\a,\b}}}$ and $x^{2^{n'_{\a,\b}+n_{s,j}+L}}=x^{2^{k_{\a,\b}}}$ for each $x \in \F$.

Now let $\tilde E_{\pm, s,j}, \tilde A_{s,j}$  be the matrices in $SL(4,R)$ obtained from their conterparts 
$E_{\pm, i_s, \pi_s(j)}^{\epsilon_s r_{s,j}}$, $A_{ i_s, \pi_s(j)}^{\epsilon_s r_{s,j}}$ by substituting each monomial from $K_{\a,\b,l}$ by the corresponding monomial in the set
\[ \tilde K_{\a,\b,l}=\left \{ t_{\a,\b,l}^{2^{b}}, \  v_{\a,\b,l}^{2^{b}} \ | \ \ b \in \{m_{\a,\b},k_{\a,\b}\} \right \}. \]
We have \begin{equation} \label{equalitytilde}\tilde A_{s,j}(\mathbf y, \mathbf z ) = \left (A_{ i_s, \pi_s(j)}(\mathbf y, \mathbf z) \right)^{\epsilon_s r_{s,j}}\end{equation} for each 
$\mathbf y=(y_1,y_2,y_3)$, $\mathbf z=(z_1,z_2,z_3) \in \F^{(3)}$.

Now define \[ \tilde U_j:= \prod_{s=1}^k \tilde A_{s,j} \in SL(4,R).\]

Every $\tilde E_{\pm, s,j,l}$ is matrix with entries which can be written as polynomials from $\F[\tilde K_{\a,\b,l}]$ of degree at most 2 in each variable from $\tilde K_{\a,\b,l}$. As $|C_{\a,\b}| \leq k$ we deduce that the degree of any coefficient of $\tilde U_j$ in any variable $t_{\a,\b,l}$ or $v_{\a,\b,l}$ is at most 
\[ 2k \max\{ 2^{k_{\a,\b}}, 2^{m_{\a,\b}}\} \leq 2k2^{2L-1-2k} =2k \cdot 4^{-k}|\F| < |\F|. \] At the same time since the elements $\tilde A_{s,j}$ satisfy the conditions from Corollary \ref{suzapplic} we deduce that $\tilde U_j$ is not central in $SL(4,R)$.
The argument from Case I therefore applies and we find elements $y_{\a,\b,l}, z_{\a,\b,l} \in \F$ such that $h_j^{-1}\tilde U_j(\mathbf y,\mathbf z)$ is not central in $SL(4,\F)$. We observe that (\ref{equalitytilde}) implies $U_j(\mathbf y,\mathbf z)=\tilde U_j(\mathbf y,\mathbf z)=u_j \in S$ where $\mathbf y =(y_{\a,\b,l}), \mathbf z =(z_{\a,\b,l}) \in \F^{(3dm)}$. Therefore as before we can set $a_{\a,\b}= A_{\a,\b} (\mathbf y_{\a,\b},\mathbf z_{\a,\b})$. This completes the last case and the proof of Proposition \ref{coset}. $\square$

\subsubsection*{Acknowledgement:} I am grateful to M. Liebeck and A. Shalev for bringing the anabelian groups to my attention and useful discussions and to J. Bray for saving me from error in the proof of Proposition \ref{coset} for Suzuki and Ree groups.


\begin{thebibliography}{99}

\bibitem{gow} L. Babai, N. Nikolov, L. Pyber, Product growth and mixing in finite
groups, \emph{SODA: Proceedings of the nineteenth annual ACM-SIAM symposium on Discrete algorithms San Francisco, California} (2008), 248-
257.

\bibitem{GLS} D. Gorenstein, R. Lyons, R. Solomon, \emph{The Classification of the finite simple groups Vol 3,} Mathematical Surveys and Monographs, vol. 40.3 (1997).

\bibitem{delaHarpe} P. de la Harpe \emph{Topics in Geometric group theory} Univeristy of Chicago Press, 2000.
\bibitem{LandSeitz} V. Landazuri, G. Seitz, On the minimal degrees of projective representations of the finite Chevalley groups, \emph{J. Algebra} 32, 418 - 443 (1974).
\bibitem{LS} M. Larsen, A. Shalev, and P. H. Tiep, The Waring problem for finite simple groups, \emph{Annals of Math. 174} (2011), 1885-1950.
\bibitem{levy} M. Levy, Word values in finite groups, Ph.D. thesis, Imperial College London, 2013.

\bibitem{ore} M. W. Liebeck, E. A. O'Brien, A. Shalev and P. H. Tiep, The Ore
conjecture, {\it J. Eur. Math. Soc.} {\bf 12} (2010), 939--1008.

\bibitem{MZ} C. Martinez, E. Zelmanov, Products of powers in finite simple groups \emph{Israel J. Math.} 96 Part B (1996) 469-479.
 
\bibitem{classicaldecomp} N. Nikolov, A product decomposition for the classical quasisimple groups,
\emph{J. Group Theory} 10 (2007), 43-53.

\bibitem{DanNN} N. Nikolov, D. Segal, On finitely generated profinite groups, I: strong completeness and uniform bounds, \emph{Annals of Math.}, {\bf 165}, (2007), 171--238.
\bibitem{NS} N. Nikolov, D. Segal, Powers in finite groups, \emph{Groups Geom. Dynamics}, Vol. 5, 501-507, 2011.
\bibitem{DanN12} N. Nikolov, D. Segal, Generators and commutators in finite groups; abstract quotients of compact groups, \emph{Invent. Math.} {\bf 190}, (2012), pp. 513-602.

\bibitem{SW} J. Saxl, J. Wilson, A note on powers in finite simple groups, \emph{Math. Proc. Camb. Phil. Soc} (1997), 91-94.
\bibitem{dan} D. Segal, \emph{Words: notes on verbal width in groups}, London Math. Soc.
Lecture Notes Series \textbf{361}, Cambridge Univ. Press, 2009.

\bibitem{shalev} A. Shalev: Words and Groups, in L. Lov\'{a}sz, I. Ruzsa, V.T. Sos, D. Palvolgyi (Eds.),\ \emph{Erd\"{o}s Centennial}, Bolyai Society Mathematical Studies, Vol. 25
2013.
\end{thebibliography}
\end{document}